\documentclass{birkjour}
\usepackage{amscd}
\usepackage{amssymb}
\usepackage{mathrsfs}
 \newtheorem{thm}{Theorem}[section]
 \newtheorem{cor}[thm]{Corollary}

 \theoremstyle{definition}
 
 \theoremstyle{remark}

 \numberwithin{equation}{section}
\usepackage[final]{hyperref}
\hypersetup{unicode= false, colorlinks=true, linkcolor=blue,
anchorcolor=blue, citecolor=red, filecolor=red, menucolor=blue,
pagecolor=blue, urlcolor=blue}
\usepackage[none]{hyphenat}
\renewcommand{\theenumi}{\roman{enumi}}
\renewcommand{\labelenumi}{(\theenumi)}

\newcommand{\CC}{\mathbb {C}}

\title[On  Toeplitz operators between Fock spaces]
 { On Toeplitz operators between  Fock spaces}
\author [Tesfa  Mengestie]{Tesfa  Mengestie }
\address{Department of Mathematical Sciences\\
Norwegian University of Science and Technology (NTNU)\\
 NO- 7491 Trondheim, Norway}
\email{tesfantnu@gmail.com}
\begin{document}
\begin{abstract}
We study mapping properties of Toeplitz operators $T_\mu$ associated to nonnegative Borel
measure $\mu$ on the complex space $\CC^n$.  We, in particular, describe  the bounded and compact operators  $T_\mu$
acting   between Fock spaces   in terms of
the objects   $t$-Berezin transforms, averaging functions, and
averaging sequences of their  inducing measures $\mu$. An  asymptotic  estimate
for the norms of  the operators   has been also obtained. The results obtained extend a recent
work of Z. Hu and X. Lv and fills the remaining gap when both the smallest and largest  Banach--Fock spaces are taken into account.
\end{abstract}
 \maketitle
\section{Introduction}
For a positive parameter  $ \alpha, $ the classical  Fock space $\mathcal{F}_\alpha^p$ consists of entire functions $f$ for which
\begin{align*}
\|f\|_{(p,\alpha)}^p=  \Big(\frac{\alpha p}{2\pi}\Big)^n\int_{\CC^n}
|f(z)|^p e^{-\frac{\alpha p}{2}|z|^2} dV(z) <\infty
\end{align*} where $ 1\leq p <\infty$ and $dV(z)$  denotes the usual Lebesgue measure
on $\CC^n.$  For $p= \infty,$ the corresponding  growth type Fock space contains those entire functions  $f$ for which
\begin{align}
\label{def}
\|f\|_{(\infty,\alpha)}= \sup_{z\in \CC^n}
|f(z)|e^{-\frac{\alpha}{2}|z|^2} <\infty.
\end{align}
 The space $\mathcal{F}_\alpha^2$ is a reproducing kernel Hilbert
space with kernel and normalized kernel functions respectively
$K_{w}(z)= e^{\alpha \langle z, w\rangle}$ and
$k_{w}(z)= e^{\alpha\langle z, w\rangle-\alpha|w|^2/2}$ where
\begin{align*}\langle z, w\rangle= \sum_{j= 1}^n z_j\overline{w_j},\ \  |z|= \sqrt{\langle z, z\rangle}, \text{and} \ w=(w_j), z=(z_j) \in \CC^n.  \end{align*}
For each $\delta>0$, we denote by $D(z, \delta)$ the  ball $\{w\in
\CC^n: |z-w|<\delta\},$ and consider a sequence of points $z_k$ such
that for $r>0,$ the sequences of the  balls  $ D(z_k, r/2)$ covers the entire space
$\CC^n$ while the  balls $ D(z_k, r/4)$ are mutually disjoint. Because of this covering property, it has been known that for any $\delta>0,$  there
exists a positive integer $N_{\max}= N_{\max}(\delta, r)$ such that every point in
$\CC^n$ belongs to at most $N_{\max}$ of the balls  $D(z_k, \delta)$. The
sequence $z_k$ will refer to such  fixed $r/2$-lattice throughout the
paper.\\
For a nonnegative measure $\mu$ on $\CC^n,$ its average  on $D(z,r)$ is the quantity
$\mu(D(z,r))/|D(z,r)|$ where $|D(z,r)|$ refers to the volume of the
ball $D(z,r)$ which depends only on $r$ for all $z$ in $\CC^n.$ Thus, we call
$\mu(D(z,r))$ simply an averaging function of $\mu,$ and $\mu(D(z_k,
r))$ its averaging sequence.

Recent years have seen a lot of  work on  Toeplitz operators acting on different
spaces of holomorphic functions since introduced by O. Toeplitz \cite{OT} in the
year 1911.  The characterization
of symbols $\psi$ in $ L^\infty$ which induce bounded or compact
Toeplitz operators $T_\psi$ on the classical Fock space $\mathcal{F}_\alpha^2$
have been studied by several authors; for example see \cite{BA1, BA2, BA3,RDV, SS,
SK}. In \cite{JIKZ},  J. Isralowitz and K. Zhu took the case further and  studied  conditions on positive Borel  measure $\mu$ on the complex plane $\CC$  under which the  induced map  $T_\mu$,
\begin{align}
\label{teopdef} T_\mu f(z)= \int_{\CC} K_{w}(z)f(w)
e^{-\alpha |w|^2} d\mu(w),
\end{align} becomes  bounded and compact on $\mathcal{F}_\alpha^2$ whenever  the measure $\mu$ satisfies  an  admissibility condition
 \begin{align}
  \label{admiss}
  \int_{\CC} |K_w(z)|^2 e^{-\alpha |w|^2} d\mu(w) <\infty.
 \end{align}
  Inspired by  this  work,  Z. Hu and
X. Lv \cite{ZHXL} recently  identified the bounded and compact
$T_\mu$  when it acts  between the Fock spaces $
\mathcal{F}_\alpha^p$ and $ \mathcal{F}_\alpha^q$ in terms of Fock--Carleson measures  whenever the exponents
$p$ and $q$ are limited  in the range $1<p,q<\infty$ and $\mu$ still satisfies  condition \eqref{admiss} on $\CC^n$. We
mention that this condition along with density of the kernel functions or the atomic decomposition result \cite[Theorem~8.2]{SJR} and Cauchy--Schwarz inequality  ensure that $T_\mu$ is densely defined on Fock spaces.

The main purpose of this note is to  extend the  results in \cite{ZHXL} for the two extreme exponents $1$ and $\infty$ assuming the   same
admissibility condition \eqref{admiss} with $\CC$ replaced by $\CC^n$ throughout the paper. It may be noted that the corresponding Fock spaces  $\mathcal{F}_\alpha^1$ and  $\mathcal{F}_\alpha^\infty$ are respectively the smallest and largest Banach--Fock spaces in  the sense of  inclusion \cite[Theorem 7.2]{SJR}.  The proofs of some of  Hu and Lv  results rely on  duality and complex interpolation arguments in Fock spaces. Because of the complicated nature  and size of the dual space of $\mathcal{F}_\alpha^\infty$, carrying over those arguments to these extreme exponent cases seems more difficult especially when the growth space $\mathcal{F}_\alpha^\infty$ is mapped into
smaller spaces $F_\alpha^p$. Thus, we will  take here  a different  approach in dealing with these two cases.

 The notion of Carleson measures has  been one of the  basic tools in studying properties of several operators  on spaces of holomorphic functions. Since the notion plays a role in our case too,  we close this section  by  recalling a few words about it.  For  $0<p\leq \infty$ and $0<q<\infty$, we call  a nonnegative measure $\mu$ on $\CC^n$ a
$(p, q) $ Fock--Carleson measure if the inequality \footnote{A word on notation:  The notation $U(z)\lesssim
V(z)$ (or equivalently $V(z)\gtrsim U(z)$) means that there is a
constant $C$ such that $U(z)\leq CV(z)$ holds for all $z$ in the set
in question, which may be a Hilbert space or  a set of complex
numbers. We write $U(z)\simeq V(z)$ if both $U(z)\lesssim V(z)$ and
$V(z)\lesssim U(z)$. We also write $L^p$ to mean $L^p(\CC^n, dV)$.}
\begin{equation}
\label{carleson} \int_{\CC^n} |f(z)|^q e^{\frac{-\alpha
q}{2}|z|^2} d\mu(z)\lesssim \| f\|_{(p,\alpha)}^q
\end{equation}holds  for all $f$ in $\mathcal{F}_\alpha^p$.  Thus $\mu$ is a (p,q) Fock--Carleson measure if and only if the
embedding map $I_\mu:\mathcal{F}_\alpha^p \to L^q(\sigma_q)$ is bounded where
$d\sigma_q(z)= e^{-q\alpha|z|^2/2} d\mu(z).$ We
 call $\mu$ a $(p, q) $ vanishing Fock--Carleson measure if such an  embedding map is
 compact.   For finite exponents $p$ and $q$   such measures were described in \cite{ZHXL}. On the other hand, the
 $(\infty, q)$ Fock--Carleson measures  were  recently identified by the author in  \cite{TM}. The description of such measure $\mu$  involves its $t$-Berezin transform  defined by
  \begin{equation*}
\widetilde{\mu}_{t}(w)= \int_{\CC^n}
|k_{w}(z)|^{t} e^{-\frac{t\alpha}{2}|z|^2}d\mu(z), \ \  t>0.
\end{equation*}
 \section{Main results}
In this section we  formulate  our results on  bounded and compact Toeplitz operators $T_\mu$
in terms of Fock--Carleson measures.  Then the corollaries  to follow give
comprehensive lists of the different characterizations of these
maps in terms of t-Berezin transforms, averaging functions, and averaging sequences of the symbols inducing
the maps.
\begin{thm} \label{bounded1}
Let $1\leq p<\infty$ and $\mu$ be a nonnegative measure on $\CC^n$.  Then $T_\mu: \mathcal{F}_\alpha^p\to
\mathcal{F}_\alpha^\infty $ is
\begin{enumerate}
\item  bounded if and only $\mu$ is a $(1,q)$
Fock--Carleson measure for some (or any) finite $q\geq p$.
Furthermore,  we have
\begin{equation} \label{normmm}
 \|T_\mu\| \simeq \|I_\mu\|
\end{equation} where $I_\mu$ is the embedding map $I_\mu:\mathcal{F}_\alpha^1 \to L^1(\sigma_1).$
\item  compact if and only if $\mu$ is a $(1,q)$
vanishing Fock--Carleson measure for some (or any) finite $q\geq p$.
\end{enumerate}
\end{thm}
While Theorem~\ref{bounded1} gives conditions under which $T_\mu$ maps smaller spaces into the  largest space $\mathcal{F}_\alpha^\infty$,  the following  result gives a stronger condition on $\mu$ when it maps  $\mathcal{F}_\alpha^\infty$ into smaller Fock spaces $\mathcal{F}_\alpha^p$ for finite $p$.
\begin{thm}\label{bounded2}
Let $1\leq p<\infty$ and $\mu$ be a nonnegative measure on $\CC^n$.  Then the following statements
are equivalent.
\begin{enumerate}
\item
 $T_\mu: \mathcal{F}_\alpha^\infty \to \mathcal{F}_\alpha^p $ is bounded;
 \item $T_\mu: \mathcal{F}_\alpha^\infty \to \mathcal{F}_\alpha^p $ is compact;
\item  $\mu$ is an $(\infty,q)$ vanishing
Fock--Carleson measure for some (or any)  finite $q\geq 1$.
 Furthermore, we have
\begin{align}
\label{normest1} \|T_\mu\| \simeq \|I_\mu\|
\end{align} where $I_\mu$ is the embedding map $I_\mu:\mathcal{F}_\alpha^\infty \to L^p(\sigma_p).$
 \end{enumerate}
\end{thm}
The phenomena with condition (iii) is, in particular,  interesting. It  ensures that if there exists a finite $q\geq1$ for which
$T_\mu: \mathcal{F}_\alpha^\infty \to \mathcal{F}_\alpha^q $ is bounded (compact), then the same conclusion holds for  $T_\mu$  when we replace the target space by $\mathcal{F}_\alpha^p$ for  all finite $p\geq 1.$ This phenomena does not necessarily happen when $T_\mu$ maps $\mathcal{F}_\alpha^p$ into $\mathcal{F}_\alpha^q$ for finite exponents $p$ and $q$  even if the target exponent $q$ is smaller than  $p.$

As consequences of the characterizations of the Fock--Carleson
measures in  Theorems~3.1  and 3.2 of \cite{ZHXL}, and
from a particular case of Theorem~2.4 of \cite{TM},  we obtain several other
characterizations for the bounded and compact operators $T_\mu$ in
terms of t-Berezin transform, averaging functions, and averaging
sequences of $\mu$. For readers convenient, we populate these conditions in the
next two corollaries.
\begin{cor} \label{cor1}
Let $1\leq p<\infty$ and $\mu$ be a nonnegative measure on $\CC^n$.  Then the following statements
are equivalent.
\begin{enumerate}
\item
 $T_\mu: \mathcal{F}_\alpha^p\to
\mathcal{F}_\alpha^\infty $ is bounded;
\item  $\mu$ is a $(1,q)$
Fock--Carleson measure for some (or any) finite $q\geq p$;
\item $\mu(D(., \delta))\in L^\infty$ for some (or any )
$\delta >0;$
\item $\widetilde{\mu_t}\in L^\infty$ for some (or any )
$t >0;$
\item $\mu(D(z_k, r))\in \ell^\infty$ for some (or any )
$r >0.$ Furthermore, we have
\begin{align}
\label{normest1} \|T_\mu\| \simeq \|I_\mu\| \simeq
\|\widetilde{\mu_t}\|_{L^\infty}\simeq \|\mu
(D(.,\delta))\|_{L^\infty}\simeq \|\mu(D(z_k,
r))\|_{\ell^\infty}
\end{align} where $I_\mu$ is the embedding map $I_\mu:\mathcal{F}_\alpha^1 \to L^1(\sigma_1).$
 \end{enumerate}
\end{cor}
The last four conditions in the corollary  above are proved in
\cite{ZHXL} and \cite{TM} where in the latter a more general setting is considered. The last three  conditions are  in addition
independent of the exponents $p$ and $\infty$ under consideration,
and remain valid if we replace the target space  $\mathcal{F}_\alpha^\infty $ by a smaller
space $\mathcal{F}_\alpha^q$ for any $q\geq p.$ Thus, if need be any one of these
conditions might be used to define the $(\infty, \infty)$
(vanishing) Fock--Carleson measures since an inequality of
the form in \eqref{carleson} makes no sense as it stands. Taking the "little oh"  counter part of the above conditions gives a list of
different characterizations of the compact Toeplitz operator $T_\mu$ acting between $\mathcal{F}_\alpha^p$
and $ \mathcal{F}_\alpha^\infty $.
\begin{cor} \label{cor3}
Let $1\leq p<\infty$ and $\mu$ be a nonnegative measure on $\CC^n$.  Then the following statements
are equivalent.
\begin{enumerate}
\item
 $T_\mu: \mathcal{F}_\alpha^\infty \to \mathcal{F}_\alpha^p $ is bounded;
 \item $T_\mu: \mathcal{F}_\alpha^\infty \to \mathcal{F}_\alpha^p $ is compact;
\item  $\mu$ is an $(\infty,q)$
Fock--Carleson measure for some (or any) finite $q\geq 1$;
\item  $\mu$ is an $(\infty,q)$ vanishing
Fock--Carleson measure for some (or any) finite $q\geq 1$;
\item $\widetilde{\mu_t}\in L^1$ for some (or any )
$t >0;$
\item $\mu$ is a finite measure on $\CC^n;$
\item $T_\mu: \mathcal{F}_\alpha^2\to \mathcal{F}_\alpha^2$ belongs to the trace class;
\item $\mu(D(., \delta))\in L^1$ for some (or any )
$\delta >0;$
\item $\mu(D(z_k, r))\in \ell^1$ for some (or any )
$r >0.$ Furthermore, we have
\begin{align}
\label{normest1}
 \|T_\mu\|^p \simeq \|I_\mu\|^p\simeq
\|\widetilde{\mu_t}\|_{L^1}\simeq \|\mu
(D(.,\delta))\|_{L^1}\simeq \|\mu(D(z_k,
r))\|_{\ell^1} \simeq \mu(\CC^n)
\end{align} where $I_\mu$ is the embedding map $I_\mu:\mathcal{F}_\alpha^\infty \to L^p(\sigma_p).$
 \end{enumerate}
\end{cor}
We note in passing that the equivalence of conditions (iii), (vi), (v), (vi), (viii) and
(ix) has been proved in \cite{TM}. On the other hand the equivalencies of (vi) and (vii) comes from
Proposition~10 in \cite{JIKZ}. These two  conditions follow from condition (v) when
we in particular set $t=1.$

We now return to the case when the exponent $p=1$ and study properties of
$T_\mu$ acting between $\mathcal{F}_\alpha^1$ and $ \mathcal{F}_\alpha^p $ and conversely.
\begin{thm}\label{particluar}
Let  $\mu$ be a nonnegative measure on $\CC^n$. Then
\begin{enumerate}
 \item  $T_\mu: \mathcal{F}_\alpha^1 \to \mathcal{F}_\alpha^p, \ 1\leq p\leq \infty, $ is \\
 a) bounded if and only if
  $\mu$ is a $(1, q)$ Fock--Carleson measure for some (or any) $q\geq p.$
  Furthermore, we have
 \begin{align}
 \left\|T_\mu\right\| \simeq \left\|I_\mu\right\|
 \end{align} where $I_\mu$ is the embedding map $I_\mu:\mathcal{F}_\alpha^1 \to L^p(\sigma_p).$\\
 b) compact if and only if $\mu$ is a $(1, q)$ vanishing Fock--Carleson measure  for some (or any) $q\geq p.$
 \item $T_\mu: \mathcal{F}_\alpha^p \to \mathcal{F}_\alpha^1,\ 1<p\leq \infty,  $ is bounded (compact) if and only if $\mu$ is a $(1,p)$ Fock--Carleson measure. We also have the norm estimate
  \begin{align}
  \label{forone}
 \left\|T_\mu\right\| \simeq \left\|I_\mu\right\|
 \end{align} where $I_\mu$ is the embedding map $I_\mu:\mathcal{F}_\alpha^p \to L^1(\sigma_1).$
 \end{enumerate}
\end{thm}
As in the previous corollaries, the conditions in Theorem~\ref{particluar} can be also  equivalently expressed in terms of  the  $t$- Berezin transform of $\mu$, the averaging functions $\mu(D(., r))$,  and  averaging sequences $\mu(D(z_j,r))$. The reformulation in terms of these notions  again follows  easily from  the characterizations of the Fock--Carleson measures in  Theorems~3.1  and 3.2 of \cite{ZHXL}, and Theorem~2.4 of \cite{TM} as  refereed  above.
 \section{Proof of the results}
  \emph{Proof of Theorem~\ref{bounded1}.} (i) We may first assume that $T_\mu: \mathcal{F}_\alpha^p
\to \mathcal{F}_\alpha^\infty$ is bounded and proceed to show the necessity of the condition. By the reproducing property of the
kernel, we have
\begin{align}
\label{kernel} K_{w}(z) =\Big(\frac{\alpha }{\pi}\Big)^n \int_{\CC^n}K_{w}(\eta)
K_{\eta}(z) e^{-\alpha|\eta|^2}dV(\eta).
\end{align} By condition \eqref{admiss}, we
observe that   the 2-Berezin transform
\begin{align*}
\widetilde{\mu}_2(z)&=\int_{\CC^n} |K_{w}(z)|^2
e^{-\alpha(|z|^2+ |w|^2)}d\mu(w)\\
&= e^{-\alpha|z|^2} \int_{\CC^n} |K_{w}(z)|^2
e^{-|w|^2)}d\mu(w)\lesssim e^{-\alpha|z|^2}  <\infty.
\end{align*}
On the other hand, our assumption ensures
\begin{align*}
\bigg|\int_{\CC^n} T_\mu
k_{z}(\eta)\overline{k_{z}(\eta)}
e^{-\alpha|\eta|^2}dV(\eta)\bigg| &= \bigg|\int_{\CC^n} \bigg(T_\mu
k_{z}(\eta) e^{-\frac{\alpha}{2}|\eta|^2}\bigg)\overline{k_{z}(\eta)}
e^{-\frac{\alpha}{2}|\eta|^2}dV(\eta)\bigg|\\
&\leq \|T_\mu
k_{z}\|_{(\infty, \alpha)} \int_{\CC^n} e^{-\frac{\alpha}{2}|\eta|^2}dV(\eta)\\
&\lesssim \int_{\CC^n}|k_{z}(\eta)| e^{-\frac{\alpha}{2}|\eta|^2}dV(\eta)\\
 &\simeq \|k_z\|_{(1,\alpha)}= 1 <\infty
\end{align*} for all $z$ in $\CC^n.$\\
Now, we may write the 2-Berezin transform as
\begin{align*}
\widetilde{\mu}_2(z)=\int_{\CC^n} |K_{w}(z)|^2
e^{-\alpha(|z|^2+ |w|^2)}d\mu(w) \ \ \ \ \ \ \ \ \ \ \ \  \ \ \ \ \
\ \ \ \ \ \ \ \ \ \ \ \ \ \ \ \ \ \ \ \ \ \ \ \ \ \ \ \ \ \ \
\end{align*}
\begin{align*}
=\int_{\CC^n} K_{w}(z) e^{-\frac{\alpha}{2}(|z|^2+ |w|^2)}
K_{z}(w) e^{-\frac{\alpha}{2}(|z|^2+ |w|^2)}d\mu(w).
\end{align*}
 Applying the kernel property in \eqref{kernel}, we
find that  the  integral above is equal to
\begin{align*}
\int_{\CC^n} K_{z}(w) e^{-\frac{\alpha}{2}(|z|^2+
|w|^2)}\bigg(\Big(\frac{\alpha }{\pi}\Big)^n\int_{\CC^n}\frac{K_{w}(\eta) K_{\eta}(z)}{
e^{\alpha|\eta|^2}} dV(\eta)\bigg)e^{-\frac{\alpha}{2}(|z|^2+
|w|^2)}d\mu(w)
\end{align*}
\begin{align*}
 =\Big(\frac{\alpha }{\pi}\Big)^n\int_{\CC^n}\int_{\CC^n}\bigg(\frac{K_{z}(w)
K_{w}(\eta)}{
e^{\alpha|w|^2+\alpha|z|^2}}\bigg)d\mu(w)\overline{K_{z}(\eta)}e^{-\alpha|\eta|^2}dV(\eta)
\end{align*}
\begin{align}
= \Big(\frac{\alpha }{\pi}\Big)^n\int_{\CC^n} T_\mu
k_{z}(\eta)\overline{k_{z}(\eta)}
e^{-\alpha|\eta|^2}dV(\eta)
\label{GG}.
\end{align} Note that for each $f$ in $ \mathcal{F}_\alpha^\infty,$   \eqref{def} implies
\begin{align}
\label{holder} |f(z)|\leq  \|f\|_{(\infty,
\alpha)}e^{\frac{\alpha}{2} |z|^2}\end{align} from which it follows
\begin{align*}
 \widetilde{\mu}_2(z) \simeq \int_{\CC^n} T_\mu
k_{z}(\eta)\overline{k_{z}(\eta)}
e^{-\alpha|\eta|^2}dV(\eta)\ \ \ \ \ \ \ \ \ \ \ \  \ \ \ \ \ \ \ \
\ \ \ \ \ \ \ \ \ \ \ \ \ \ \ \ \ \ \ \ \ \ \ \ \ \ \ \ \ \ \ \ \ \
\end{align*}
\begin{align}\label{equall}
\leq \|T_\mu
k_{z}\|_{(\infty,\alpha)}\int_{\CC^n}|k_{z}(\eta)|e^{-\frac{\alpha}{2}
|\eta|^2}dV(\eta) \lesssim \|T_\mu\|,
\end{align}
 where  the
last inequality follows because  $\|k_z\|_{(1,\alpha)}=\|k_z\|_{(\infty,\alpha)}=1$ for all
$z$ in $\CC^n.$  This proves that $\widetilde{\mu_2}$
is uniformly  bounded on $\CC^n.$   By
 Theorem 3.1 of \cite{ZHXL}, this  happens  if and only if $\mu$ is a  $(p,q)$ Fock--Carleson measure
 for any   $p$ and $q$ in the range $0<p\leq q<\infty,$ from which the desired conclusion  follows. In particular when  $p=q=1$ the result
 there  and \eqref{equall} imply
 \begin{align}
 \label{othersidee}
  \sup_{z\in
\CC^n}\widetilde{\mu_2}(z)\simeq \|I_\mu\| \lesssim \|T_\mu\|
\end{align} where $I_\mu$ is the embedding map $I_\mu:\mathcal{F}_\alpha^1 \to L^1(\sigma_1).$\\
On the other hand, if  $\mu$ is a  $(1,q)$ Fock--Carleson measure
for some finite  $q\geq p\geq1$, then it is a $(1,1)$ Fock--Carleson measure; see Theorem~3.1 of \cite{ZHXL} again. Thus,
\begin{align}
\label{normof}
 |T_\mu f(z)| e^{-\frac{\alpha}{2}|z|^2}&\leq
\int_{\CC^n} |k_{z}(w)| |f(w)|e^{-\alpha
|w|^2}d\mu(w)\nonumber\\
&\leq \|f\|_{(\infty,\alpha)} \int_{\CC^n}
|k_{z}(w)|e^{-\frac{\alpha}{2}|w|^2} d\mu(w)\\
&\lesssim \left\|I_\mu\right\|\|f\|_{(\infty,\alpha)}\|k_{z}\|_{(1,\alpha)}=\left\|I_\mu\right\|\|f\|_{(\infty,\alpha)}.\label{OO}
\end{align}Here  the second inequality follows  by \eqref{holder}. From this it follows that
\begin{align*}
\sup_{z\in \CC^n} |T_\mu f(z)| e^{-\frac{\alpha}{2}|z|^2} \lesssim
\left\|I_\mu\right\| \|f\|_{(\infty,\alpha)}\lesssim \left\|I_\mu\right\|\|f\|_{(p,\alpha)}
\end{align*}where we use the inclusion $\mathcal{F}_\alpha^p \subseteq
\mathcal{F}_\alpha^\infty,$ and hence $T_\mu$ is bounded and
\begin{align}
\label{pp}
\|T_\mu\| \lesssim \left\|I_\mu\right\|.
\end{align}
From   \eqref{othersidee} and \eqref{pp},  the asymptotic norm estimate in
\eqref{normmm} follows.

 (ii) We may first verify the sufficiency of the condition.
   Let $f_m$ be a sequence of
 functions in $\mathcal{F}_\alpha^p$ such that
$\sup_m \|f_m\|_{(p,\alpha)}<\infty$ and $f_m$ converges uniformly
to zero on compact subsets of $\CC^n$ as $m\to \infty$. An application of   Lemma~1 of \cite{SSS} gives
 \begin{align*}
|T_\mu f_m(z)| e^{-\frac{\alpha}{2}|z|^2}\leq\sup_{m\geq1}
\|f_m\|_{(p,\alpha)} \int_{\CC^n }|k_{z}(w)|
e^{-\frac{\alpha}{2} |w|^2}d\mu(w)\\
\lesssim  \int_{\CC^n }|k_{z}(w)| e^{-\frac{\alpha}{2}
|w|^2}d\mu(w).
\end{align*}
Since $\mu$ is a $(1,q)$ vanishing Fock--Carleson measure for some finite $ q\geq p\geq1$, by Theorem~3.2 of \cite{ZHXL},  it is also a $(1,1)$ vanishing Fock--Carleson measure.  This along with the fact that  $k_{z}$
converges uniformly to zero on compact subsets of $\CC^n$ when $|z|
\to \infty$   yields
\begin{align*}\int_{\CC^n }|k_{z}(w)| e^{-\alpha
|w|^2/2}d\mu(w) \to 0\end{align*} as $|z| \to \infty$ from which we conclude
that $T_\mu$ is compact.

 Conversely, suppose that
$T_\mu:\mathcal{F}_\alpha^p \to \mathcal{F}_\alpha^\infty$ is compact. By arguing  with \eqref{GG} and  \eqref{holder}, we have
\begin{align*} \widetilde{\mu_2}(z)=\int_{\CC^n} |k_{z}(w)|^2 e^{- \alpha|w|^2}
d\mu(w)
\leq \|T_\mu k_{z}\|_{(\infty,\alpha)}\|
k_{z}\|_{(1,\alpha)}\\
=\|T_\mu k_{z}\|_{(\infty,\alpha)}
 \to 0
\end{align*} as $|z| \to \infty$ since  $k_{z}$ is a sequence of
unit norm functions which converges uniformly to zero on compact subsets of $\CC^n$ when $|z|\to \infty$. Then, by a particular case of  Theorem~3.2 of
\cite{ZHXL}, it follows that $\mu$ is a $(1,q)$ vanishing
Fock--Carleson measure for any finite $q\geq p$.

\emph{Proof of Theorem~\ref{bounded2}.}  Since (ii) implies
(i) is trivial, we only need to  verify (iii) implies (ii)  and
(i) implies (iii). Suppose $\mu$ is an $(\infty, q)$
vanishing Fock--Carleson measure for some
 finite $q\geq 1$. Then for each $f$ in $\mathcal{F}_\alpha^\infty$,
 \begin{align}
 \label{twopart}
 \left\|T_\mu f\right\|_{(p,\alpha)}^p &\simeq \int_{\CC^n} \left|T_\mu f(z)e^{-\frac{\alpha}{2}|z|^2} \right|^pdV(z)\nonumber\\
 &\leq \int_{\CC^n} \bigg( \int_{\CC^n}\frac{ |K_w(z)| |f(w)|}{e^{\alpha|w|^2+\frac{\alpha}{2}|z|^2}}  d\mu(w)\bigg)^p dV(z)= I.
  \end{align}
 We may now first assume that $p>1.$ Then  applying  H\"{o}lder's inequality on the right-hand side of \eqref{twopart} and subsequently Fubini's theorem give
  \begin{align}
  \label{vanishing}
 I &\leq\mu(\CC^n)^{\frac{p}{p'}} \int_{\CC^n} \left|f(w)e^{-\frac{\alpha}{2}|w|^2}\right|^p \int_{\CC^n}  |K_w(z)|^p e^{-\frac{p\alpha}{2}(|w|^2+ |z|^2)} dV(z)d\mu(w)\nonumber\\
 &\simeq  \mu(\CC^n)^{\frac{p}{p'}} \int_{\CC^n} \left|f(w)e^{-\frac{\alpha}{2}|w|^2}\right|^p d\mu(w)
 \end{align} where  $p$ and $p'$ are conjugate exponents, and by the fact that
  $\left\|k_{w} \right\|_{(p,\alpha)}= 1$ for each $p$ and $w$ in $ \CC^n$.
 On the other hand  for $p=1$, it follows that
 \begin{align}
 \label{partly}
 I=\int_{\CC^n}  \int_{\CC^n} |K_w(z)|e^{-\alpha|w|^2-\frac{\alpha}{2}|z|^2}  |f(w)|d\mu(w) dV(z)\nonumber\\
 \simeq \int_{\CC^n}    |f(w)| e^{-\frac{\alpha}{2}|w|^2}d\mu(w)
 \end{align} after interchanging the integrals.
 Since $\mu$ is assumed to be  an $(\infty,q)$ vanishing Fock--Carleson measure, a particular case of Theorem~2.4 of \cite{TM} ensures that $\mu(\CC^n)$ is finite. To this end, we  combine  \eqref{twopart}, \eqref{vanishing} and \eqref{partly} to conclude the norm estimate
  \begin{align}
\label{import}
 \|T_\mu f\|_{(p,\alpha)}^p  \lesssim
\int_{\CC^n} |f(z) e^{-\frac{\alpha}{2} |z|^2}|^p d\mu(z).
 \end{align}
 We may now
consider  a sequence of functions $ f_m$ in $\mathcal{F}_\alpha^\infty$ such
that
 $\sup_m \|f_m\| _{(\infty, \alpha)}<\infty$ and $f_m$ uniformly
 converges to zero on compact subsets of $\CC^n$ as $m \to \infty.$
 Applying the estimate in \eqref{import}, we have
 \begin{align*}
 \|T_\mu f_m\|_{(p,\alpha)}^p \lesssim \int_{\CC^n} |f_m(z) e^{-\alpha |z|^2/2}|^p d\mu(z) \to 0
\end{align*}when  $m \to \infty$. This holds  as  $\mu$ is an $(\infty, q)$ vanishing Fock--Carleson measure, by
Theorem~2.4  of \cite{TM}, it is  also an $(\infty, p)$
vanishing Fock--Carleson measure. Thus, we conclude that  $T_\mu$ is
compact. \\
Furthermore, observe that since $\mu$ is an $(\infty,p)$ Fock--Carleson measure,  inequality  \eqref{import} gives  the
estimate
\begin{align*}
\|T_\mu f\|_{(p,\alpha)}^p  \lesssim
\int_{\CC^n} |f(z) e^{-\alpha |z|^2/2}|^p d\mu(z)\leq \|I_\mu\|^p \|f\|_{(\infty,\alpha)}^p
\end{align*} for each $f$ in $\mathcal{F}_{\alpha}^\infty$ from which we find
\begin{align}
\label{on}
 \|T_\mu\| \lesssim \|I_\mu\|
 \end{align} where $I_\mu$ is the embedding map $I_\mu:\mathcal{F}_\alpha^\infty \to L^p(\sigma_p).$

It  remains to prove  (i) implies (iii).    We first consider the case $p=1$.  Suppose  $T_\mu: \mathcal{F}_\alpha^\infty
\to \mathcal{F}_\alpha^1$ is bounded. Then by Lemma~4.1 of \cite{ZHXL}
\begin{align*}
\int_{\CC^n} K_w(\eta) \overline{K_z(\eta)} e^{-\alpha|\eta|^2}d\mu(\eta)= \langle T_\mu K_w, K_z\rangle
\end{align*} for each kernel functions $K_w$ and $K_z$. From this and \cite[Theorem~8.2]{SJR}, we easily see that
for each $f$ in $\mathcal{F}_\alpha^\infty$
\begin{align}
\label{part0} \int_{\CC^n}|f(z)|^2 e^{-\alpha|z|^2} d\mu(z)
 \leq \|T_\mu f\|_{(1,\alpha)}\|f\|_{(\infty, \alpha)}
  \leq \|T_\mu\|\|f\|_{(\infty, \alpha)}^2.
\end{align} This means that  $\mu$ is an $(\infty, 2)$ Fock--Carleson measure.
Then by a particular case of Theorem~2.4 of  \cite{TM}, $\widetilde{\mu}_t\in
L^1$  and hence $\mu$ is an $(\infty, q)$ Fock--Carleson measures for  any
 $q\geq 1$ from which  the desired conclusion follows. From \eqref{part0}, we also have
 \begin{align}
 \label{neww}
 \|I_\mu\| \leq \|T_\mu\|.
 \end{align}
We now consider the case $p>1.$ For this, we may  assume that   $T_\mu: \mathcal{F}_\alpha^\infty
\to \mathcal{F}_\alpha^p$ is bounded  and   aim  to show that
 $\mu$  is a finite measure on $\CC^n.$ Then the desired conclusion  will follow
 once from a particular case of Theorem~2.4 of \cite{TM} again. To proceed further, we
 pick a suitable   test function defined by
\begin{align*}
f_0= \sum_{j=1}^\infty k_{z_j}.
\end{align*}
 Observe that
 $f_0$ belongs to $\mathcal{F}_\alpha^\infty$ since it  could be written as sum of  $c_jk_{z_j}$ where  $c_j= 1$ for all $j$ and
$\|f_0\|_{(\infty,\alpha)}\simeq 1$; see the atomic decomposition result \cite[Theorem~8.2]{SJR}. Applying $T_\mu$
 to this function gives
  \begin{align*}
 \|T_\mu\|^p&\geq \int_{\CC^n}\bigg|T_\mu \sum_{j=1}^\infty k_{z_j}(z)\bigg|^p e^{\frac{-p\alpha}{2}|z|^2} dV(z)\nonumber\\
  &=\int_{\CC^n}\Bigg|\int_{\CC^n} \sum_{j=1}^\infty e^{-\frac{\alpha|z_j|^2}{2}+\alpha\langle w,z_j\rangle + \alpha\langle z, w\rangle-\alpha|w|^2} d\mu(w)\Bigg|^p e^{\frac{-p\alpha}{2}|z|^2}dV(z) \nonumber\\
  &= \int_{\CC^n}\Bigg|\sum_{j=1}^\infty\int_{\CC^n} e^{i\alpha \Im\big(\langle w,z_j\rangle+\langle z,w\rangle\big)} e^{-\frac{\alpha|z_j-w|^2}{2}-\frac{\alpha|z-w|^2}{2} }d\mu(w)\Bigg|^p dV(z)\\
  &\simeq\int_{\CC^n}\Bigg(\sum_{j=1}^\infty\int_{\CC^n}  e^{-\frac{\alpha|z_j-w|^2}{2}-\frac{\alpha|z-w|^2}{2} }d\mu(w)\Bigg)^p dV(z)=J
   \end{align*} where $\Im \eta $  refers to  the imaginary part of a complex number $\eta$.
   Since all the terms in the last  sum are  positive and $p>1$, we also  have the estimate
   \begin{align*}
   J&\geq \int_{\CC^n}\sum_{j=1}^\infty\bigg(\int_{\CC^n}  e^{-\frac{\alpha|z_j-w|^2}{2}-\frac{\alpha|z-w|^2}{2} }d\mu(w)\bigg)^p dV(z)\nonumber\\
   &\geq \int_{\CC^n}\sum_{j=1}^\infty \int_{\CC^n} e^{-\frac{\alpha p|z_j-w|^2}{2}-\frac{\alpha p|z-w|^2}{2} }d\mu(w) dV(z)\nonumber\\
   &= \int_{\CC^n}\bigg( \int_{\CC^n} \sum_{j=1}^\infty  e^{-\frac{\alpha p|z_j-w|^2}{2}-\frac{\alpha p|z-w|^2}{2} } dV(z)\bigg)d\mu(w)\nonumber\\
   &\simeq \int_{\CC^n} \sum_{j=1}^\infty  e^{-\frac{\alpha p|z_j-w|^2}{2}}d\mu(w) \simeq \int_{\CC^n}d\mu(w)= \mu(\CC^n)
   \end{align*} where the first equality is due to Fubini's Theorem. From the above series of estimations, we deduce
   \begin{align}
  \label{last}
  \mu(\CC^n) \lesssim \|T_\mu\|^p
  \end{align} from which,  \eqref{on}, \eqref{neww} and the estimates in Theorem~2. of \cite{TM}, the asymptotic  norm relation  in
  \eqref{normest1} follows and completes the proof.

\emph{Proof of Theorem~\ref{particluar}.}
The proof of part (i) of the theorem follows from a simple  variant of  the arguments used to prove Theorem~\ref{bounded1}. Thus, we will omit it.
 To prove part (ii) of the result, it suffices  to  verify the statement: $T_\mu: \mathcal{F}_\alpha^p \to \mathcal{F}_\alpha^1$ is bounded if and only if $\mu(D(z_j,r))$ belongs to  $\ell^{p/(p-1)}$ for some positive $r$. Then the  desired conclusion  will follow once  from Theorem~3.3 of \cite{ZHXL}.  We may first assuming that  $T_\mu: \mathcal{F}_\alpha^p \to \mathcal{F}_\alpha^1$ is bounded. Note that for $p= \infty$, the result follows from Corollary~\ref{cor3} above. Thus,  we only look at the case when  $p>1$ is finite.  By  Theorem~9.3 of \cite{SJR}, our assumption  holds if and only if  $T_\mu: \mathcal{F}_\alpha^{2p} \to \mathcal{F}_\alpha^{2p-1}$ is bounded.
 Following similar arguments as those leading
 to \eqref{part0}, for each $f$ in $\mathcal{F}_\alpha^{2p}$  we  have
 \begin{align}
 \label{opp}
 \int_{\CC^n}|f(z)|^2 e^{-\alpha|z|^2} d\mu(z)
 \leq \|T_\mu f\|_{(2p-1,\alpha)}\|f\|_{(2p, \alpha)}
  \leq \|T_\mu\|\|f\|_{(2p, \alpha)}^2.
 \end{align} This proves that $\mu$ is an $(2p, 2)$ Fock--Carleson measure. This holds if and only if $\mu(D(z_j,r))$
 belongs to $\ell^{2p/(2p-2)}=\ell^{p/(p-1)} $ by special case of Theorem~2.3 of \cite{TM}.  It follows from Theorem~4.4 of \cite{ZHXL}  that  $T_\mu : \mathcal{F}_\alpha^p \to  \mathcal{F}_\alpha^1$ is bounded and this along with the series of estimations in Theorem 2.3 of \cite{TM} give
 \begin{align}
 \label{op}
 \|I_\mu\| \simeq \|\mu(D(z_j,r))\|_{\ell^{\frac{p}{p-1}}}\lesssim \|T_\mu\|
 \end{align} where $I_\mu$ is here  the embedding map $I_\mu:\mathcal{F}_\alpha^p \to L^1(\sigma_1).$
 Conversely, we assume that  $\mu(D(z_j,r))$ belongs to  $\ell^{p/(p-1)}$ for some positive $r$. Then  for each $f$ in $ \mathcal{F}_\alpha^p$, we have
   \begin{align}
   \label{suff}
  \left\|T_\mu f\right\|_{(1,\alpha)}&\leq \int_{\CC^n} \int_{\CC^n}  |f(w)| e^{-\alpha |w|^2-\frac{\alpha}{2}|z|^2}|K_{w}(z)| d\mu(w) dV(z) \nonumber\\
   &\simeq\int_{\CC^n} |f(w)| e^{-\frac{\alpha}{2} |w|^2}
    \Bigg(\int_{\CC^n} |K_w(z)| e^{-\frac{\alpha}{2}(|z|^2+|w|^2)} dV(z)\Bigg) d\mu(w) \nonumber\\
   &\simeq\int_{\CC^n} |f(w)| e^{-\frac{\alpha}{2} |w|^2} d\mu(w)
       \end{align} by Fubini's theorem and the fact that
       $\left\|k_{w} \right\|_{(p,\alpha)}= 1$ for each $p$ and $w $ in $ \CC^n$.
            By triangle inequality, it follows that
 $ D(z,r) \subset D(z_j,2r)$ whenever $z\in D(z_j,r)$. For such $z$, Lemma~1 of \cite{JIKZ} gives
 \begin{align*}
 |f(z)|e^{-\frac{\alpha}{2}|z|^2} \lesssim \int_{D(z_j, 2r)} |f(w)|e^{-\frac{\alpha}{2}|w|^2} dV(w).
 \end{align*}
  Applying this estimation to the right hand of \eqref{suff} followed by Fubini's theorem gives
  \begin{align*}
   \int_{\CC^n} |f(w)| e^{-\frac{\alpha}{2} |w|^2} d\mu(w) \lesssim
   \sum_{j=1}^\infty \mu(D(z_j,r)) \int_{D(z_j,2r)}|f(z)| e^{-\frac{\alpha}{2} |z|^2}
  dV(z).
  \end{align*}
  Applying  H\"{o}lder's inequality, we see that the right hand side quantity is bounded by
  \begin{align}
  \Bigg(\sum_{j=1}^\infty \mu(D(z_j,r))^{\frac{p}{p-1}}\Bigg)^{{\frac{p-1}{p}}}\Bigg( \sum_{j=1}^\infty \Bigg(\int_{D(z_j,2r)} |f(z)|  e^{-\frac{\alpha }{2} |z|^2} dV\bigg)^p\Bigg)^{\frac{1}{p}}\nonumber
  \end{align}
  \begin{align}
  \lesssim \left\|\mu(D(z_j,r))\right\|_{\ell^{\frac{p}{p-1}}} \Bigg(\sum_{j=1}^\infty \int_{D(z_j,2r)} |f(z)|^p  e^{-\frac{\alpha p }{2} |z|^2}dV(z)\Bigg)^{\frac{1}{p}}\nonumber
  \end{align}
  \begin{align}
  \label {last}
     \lesssim \left\|\mu(D(z_j,r/2))\right\|_{\ell^{\frac{p}{p-1}}} \left\|f\right\|_{(p,\alpha)}.
    \end{align}
     where the first   inequality follows by Jenesen's inequality as  the Lebesgue  measure is a finite measure on balls and $x^p$ is convex for all $p\geq 1.$    From \eqref{suff} and \eqref{last}, we conclude
    \begin{align}
    \left\|T_\mu\right\| \lesssim \left\|\mu(D(z_j,r/2))\right\|_{\ell^{\frac{p}{p-1}}}.
    \end{align}
    We combine this with \eqref{op} and the series of norm estimates in
    Theorem~3.3 of \cite {ZHXL} to  arrive at
    \eqref{forone} and  finish the proof.


\begin{thebibliography}{BRSHZ}
 \bibitem{BA1} Berger, C., Coburn, L.:
Toeplitz operators and quantum mechanics, J. Funct. Anal.,
\textbf{68}, 273--299 (1986)
 \bibitem{BA2}Berger, C., Coburn, L.: Toeplitz operators on the
Segal--Bargmann space, Trans. Am. Math. Soc., \textbf{301(2)}, 813--829 (1987)

\bibitem{BA3} Coburn,  L., Isralowitz,  J., Li, B.:  Toeplitz operators with
BMO symbols on the Segal--Bargmann space, Trans. Am. Math.
Soc.,\textbf{363},  3015--3030 (2011)

\bibitem{ZHXL}Hu, Z., Lv, X.: Toeplitz operators from one Fock space to
another, Integr. Equ. Oper. Theory, \textbf{70},  541--559 (2011)

\bibitem{JIKZ}Isralowitz, J., Zhu,  K.:  Toeplitz operators on the Fock
space, Integr. Equ. Oper. Theory,  \textbf{66(4)}, 593--611 (2010)

\bibitem{SJR}Janson, S., Peetre, J., Rochberg,  R.:  Hankel forms and the Fock space,  Rev. Mat. Iberoamericana,  \textbf{3}, 61--138  (1987)
\bibitem{TM} Mengestie, T.:  Carleson type measures for Fock--Sobolev spaces, to appear at Complex Analysis and Operator Theory,  DOI: 10.1007/s11785-013-0321-7
 \bibitem{RDV}Ramirez De Arellano, E., Vasilevski, E.:  Toeplitz operators on the Fock space with
pre symbols discontinuous on a thick set, Math. Nachr., \textbf{180}, 299--315 (1996)
 \bibitem{SSS}Stevi$\acute{\text{c}}$,  S.:  Weighted composition
operators on Fock-type spaces in $\CC^N$, Applied Mathematics and
Computation, \textbf{215},  2750--2760 (2009)
\bibitem{SS}Stroethoff,  S.:    Compact
Toeplitz operators on generalized Fock spaces,  Acta Sci. Math. (Szeged), \textbf{64}, 657--669   (1998)
\bibitem{SK} Stroethoff, K.:  Hankel and Toeplitz operators on the Fock space, Mich. Math. J., \textbf{39(1)},  3--16  (1992)
\bibitem{OT}Toeplitz, O.:  Zur Theorie der quadratischen und bilinearen Formen von
unendlichvielen, Ver\"{a}nderlichen. Math. Ann., \textbf{70},  351--376  (1911)
\end{thebibliography}
\end{document}